\documentclass{article}
\usepackage[utf8x]{inputenc}
\usepackage{graphicx}
\usepackage{wrapfig}
\usepackage[pdftex,bookmarks]{hyperref}
\hypersetup{
pdftitle={Boolean},
pdfauthor={José M. Zapata}
}

\usepackage{fancyhdr}
\usepackage{stmaryrd}

\usepackage{mathrsfs} 
\usepackage{relsize}

\usepackage[english]{babel}
\usepackage{amsmath}
\usepackage{amsthm}
\usepackage{amsfonts}

\usepackage{enumitem}

\usepackage{amscd}

\usepackage{epigraph}

\usepackage{yhmath}

\usepackage{amssymb}

\newcommand{\N}{\mathbb{N}}
\newcommand{\R}{\mathbb{R}}
\newcommand{\E}{\mathbb{E}}
\newcommand{\PP}{\mathbb{P}}

\newcommand{\VA}{V^{(\mathcal{A})}}
\newcommand{\VAA}{\overline{V}^{(\mathcal{A})}}
\newcommand{\auu}{\uparrow}
\newcommand{\n}{\mathfrak{n}}

\newcommand{\F}{\mathcal{F}}
\newcommand{\EE}{\mathcal{E}}
\newcommand{\A}{\mathcal{A}}

\DeclareMathOperator*{\esssup}{sup\,}

\DeclareMathOperator{\dom}{dom}
\DeclareMathOperator{\Ord}{Ord}

\newtheorem{defn}{Definition}[section]
\newtheorem{rem}{Remark}[section]
\newtheorem{prop}{Proposition}[section]

\newtheorem{thm}{Theorem}[section]

\newtheorem{example}{Example}[section]

\providecommand{\keywords}[1]{\textbf{\textit{Keywords:}} #1}

\begin{document}
\title{A Boolean valued analysis approach to conditional risk}
\author{J.~M. Zapata \thanks{Department of Mathematics, University of Konstanz, Germany. e-mail: jmzg1@um.es. } \thanks{The author was partially supported by the grants MINECO MTM2014-57838-C2-1-P and Fundaci\'{o}n S\'{e}neca 20903/PD/18.}}

\date{}
\maketitle

\begin{center}
\emph{To Evgeny Israilevich Gordon with respect and admiration on the occasion of his anniversary. }
\end{center}

\begin{abstract}
By means of the techniques of Boolean valued analysis, we provide a transfer principle between duality theory of classical convex risk measures and duality theory of conditional risk measures. 
Namely, a conditional risk measure can be interpreted as a classical convex risk measure within a suitable set-theoretic model. 
As a consequence, many properties of a conditional risk measure can be interpreted as basic properties of convex risk measures. 
 This amounts to a method to interpret a theorem of dual representation of convex risk measures as a new theorem of dual representation of conditional risk measures.  
As an instance of application, we establish a general robust representation theorem for conditional risk measures and study different particular cases of it.    

\end{abstract}

\section*{Introduction}

The present paper contributes to mathematical finance by means of the tools of Boolean valued analysis, a branch of functional analysis that applies special model-theoretic techniques to analysis. 

Let us start by explaining the mathematical finance problem that we are interested in. 
Over the past two decades and having its origins in the seminal paper \cite{artzner1999coherent}, duality theory of risk measures  has been an active and prolific area of research, see e.g.~\cite{arai2010convex,arai2014convex,biagini2009extension,xxx,cheridito02,delbaen2009risk,
farkas2014beyond,frittelli2002putting,JST06,kaina2009convex,orihuela2012lebesgue} and references therein. 
The simplest situation is the case in which only two instants of time matter: today $0$ and tomorrow $T>0$. 
In this case, the market information that will be observable at time $T$ is described by a probability space $(\Omega,\mathcal{E},\PP)$. 
A \emph{risk measure} is a function that assigns to any $\mathcal{E}$-measurable random variable $x$, which models a final payoff, a real number $\rho(x)$, which quantifies the riskiness of $x$. 
Generally speaking, duality theory of risk measures studies what are the desirable economic properties that should have a risk measure and which is the dual representation of a risk measure with these properties.

A more intricate situation is when we have a dynamic configuration of time, in which the arrival of new information at an intermediate date $0<t<T$ is taken into account.     
Suppose that the available information  at time $t$ is encoded in a sub-$\sigma$-algebra $\F$ of the $\sigma$-algebra $\mathcal{E}$ of general information. 
In that case, the riskiness at time $t$ of any final payoff  is contingent on the information contained in $\F$.  
Then a \emph{conditional risk measure} is a mapping (fulfilling some desirable economic conditions) that assigns to any final payoff, modeled by an $\mathcal{E}$-measurable random variable $x$, an $\F$-measurable random variable $\rho(x)$, which quantifies the risk arisen from $x$. 
A problem that has drawn the attention for a long time is the dual representation of a risk measure in a multi-period setup, see for instance~\cite{artzner2007coherent,biagini2014dynamic,bielecki2016dynamic,
bion2008dynamic,detlefsen2005conditional,bionnadal04,kupper11,follmer2006convex,frittelli2011conditional,frittelli14} and references therein.

As explained in \cite{kupper03}, whereas classical convex analysis perfectly applies to the one-period case, it has a rather delicate application to the multi-period model: consider the properties of the conditional
risk measure $\rho$ such as convexity, continuity, differentiability and so on. These properties have
to be satisfied by the function $x\mapsto \rho(x)(\omega)$ for each $\omega\in\Omega$, but they should be fulfilled under a measurable dependence on $\omega$, in order to enable a recursive multi-period scheme. 
This approach would require heavy measurable selection criteria.  
These difficulties have motivated some new developments in functional analysis.  
For instance, Filipovic et al.~\cite{kupper03} proposed to consider modules over $L^0(\F)$, the space of (equivalence classes of) $\F$-measurable random variables (see also~\cite{kupper11,Cheridito2012,guo10}). 
 More sophisticated machinery is provided in \cite{DJKK13}, where the so-called \emph{conditional set theory} is introduced and developed.    
Other related approaches are introduced in \cite{kuppermaccheroni,CerreiaVioglio2017,eisele13}.

A step forward is given in \cite{avile2017boolean}, where it is established a method to interpret any theorem of convex analysis as a theorem of $L^0$-convex analysis. 
The machinery is taken from Boolean valued analysis, a branch of functional analysis that  consists in studying  the properties of a mathematical object by interpreting it as a simpler object in a different set-theoretic model whose construction utilizes
a Boolean algebra. 
 Boolean valued analysis stems from the method of forcing that Paul Cohen created to prove the independence of the continuum hypothesis from the system of axioms of the Zermelo-Fraenkel set
theory with the Axiom of Choice (ZFC) \cite{cohen1966set}. 
 The main tool of Boolean valued analysis are Boolean valued models of set theory, which were developed by Scott, Solovay, and Vop\v{e}nka as a way to simplify the Cohen's method of forcing. 
 Boolean valued analysis started with Gordon~\cite{gordonReals} and Takeuti~\cite{takeuti1978}\footnote{Actually the term Boolean valued analysis was coined by Takeuti~\cite{takeuti1978}.}, and has undergone a fruitful and deep development due to Kusraev and Kutateladze. For a thorough account, we refer the reader to \cite{kusraev2012boolean} and its extensive list of references. 
 
  The present paper is aimed to extend and exploit the connections provided in \cite{avile2017boolean}, to establish a general transfer method between duality theory of one-period risk measures and duality theory of conditional risk measures, putting at the disposal of mathematical finance a powerful tool to obtain different duality representation results. 
 Namely, we show that if $\rho:\mathscr{X}\rightarrow L^0(\F)$ is a conditional risk measure, then we can interpret $\rho$ as a one-period risk measure $\rho{\uparrow}$ defined on a space of (classes of equivalence of) random variables $\mathscr{X}{\auu}$ within a suitable set-theoretic model.  
  Then, inside of this model, any available theorem about the dual representation of the one-period risk measure $\rho{\auu}$ has a counterpart that is satisfied by the conditional risk measure $\rho$. 
  This means that any theorem of duality theory of one-period risk measures gives rise to a new theorem of duality theory of conditional risk measures.    

The paper is structured as follows: In Section 1, we give some preliminaries and review duality theory of risk measures both in the one-period and multi-period setups. 
In Section 2, we recall the basics of Boolean valued models. 
In Section 3, we establish a Boolean valued transfer principle between duality theory of convex risk measures and duality theory of conditional risk measures. 
By applying this transfer principle we derive a general robust representation theorem of conditional risk measures and study different particular cases.  
Finally, in Section 4, due to limited space, we sketch the proof of the transfer method.

\section{Preliminaries on duality theory of risk measures}

In this section, we review the main elements of duality theory of risk measures. 
We start by the one-period setup, recalling the notion of convex risk measure and different properties that matter in the dual representation of a convex risk measure. 
After this, we move on to the multi-period setup. We recall the notion of conditional risk measure and introduce conditional analogues of the different elements of the one-period case.

\subsection{One-period setup: convex risk measures}

Let us recall some basics of duality theory of risk measures. 
For an introduction to this topic, we refer the reader to \cite[Chapter 4]{foellmer01}.  
Let $(\Omega,\mathcal{E},\PP)$ be a probability space.   
We denote by $L^0(\mathcal{E})$ the space of $\mathcal{E}$-measurable real-valued random variables on $\Omega$ identified whenever their difference is $\PP$-negligible. 
Given $x,y\in L^0(\mathcal{E})$ we understand $x\leq y$ and $x<y$ in the almost surely sense. 
Endowed with the order $\leq$, $L^0(\mathcal{E})$ is a Dedekind complete lattice ring.   
We say that $\lim_n x_n=x$~a.s. in $L^0(\mathcal{E})$ whenever $(x_n)$ converges almost surely to $x\in L^0(\mathcal{E})$ 
(or equivalently, $x_n$ order converges to $x$). 
 
Suppose that our probability space $(\Omega,\EE,\PP)$ models the market information at some time horizon $T>0$. 
The final payoff of each financial position is going to be modeled by a subspace $\mathcal{X}$ of $L^0(\mathcal{E})$ with the following properties:

\begin{itemize}
\item $\mathcal{X}$ is \emph{solid}, that is, $y\in\mathcal{X}$ and $|x|\leq|y|$ imply that $x\in\mathcal{X}$;\footnote{Solid subspaces are also called \emph{order ideals}.}
\item $\mathbb{E}_\mathbb{P}[|x|]<\infty$ for any $x\in \mathcal{X}$;
\item  the classes of equivalence of constant functions are contained in $\mathcal{X}$.
\end{itemize}

\begin{example}
The following subspaces of $L^0(\mathcal{E})$ satisfy the properties above:
\begin{enumerate}
\item $L^p$ spaces: $L^p(\mathcal{E}):=L^p(\Omega,\mathcal{E},\PP)$ with $1\leq p\leq\infty$.
\item Orlicz spaces: Let  $\phi:[0,\infty)\rightarrow[0,\infty]$ be a \emph{Young function}, that is, an increasing left-continuous convex function finite on a neighborhood of $0$ with $\varphi(0) = 0$ and $\lim_{x\rightarrow\infty}\phi(x) = \infty$. 
The associated \emph{Orlicz space} is
\[
L^\phi(\mathcal{E})=\{x\in L^0(\EE)\colon \exists r\in (0,\infty),\:  \E_\PP[\phi(r|x|)]<\infty\}.
\] 
\item Orlicz-heart spaces: If $\phi$ is a Young function, the associated \emph{Orlicz-Heart space} is 
\[
H^\phi(\EE)=\{x\in L^0(\EE)\colon \forall r\in (0,\infty),\: \E_\PP[\phi(r|x|)]<\infty\}.
\]  
\end{enumerate}

\end{example}

The riskiness of any final payoff $x\in\mathcal{X}$ is quantified by a function $\rho:\mathcal{X}\rightarrow\R$ satisfying for all $x,y\in\mathcal{X}$:
\begin{enumerate}
\item\emph{Convexity}: i.e. $\rho(r x + (1-r)y)\leq r\rho(x)+(1-r)\rho(y)$ for all $r\in\R$ with $0\leq r\leq 1$;
\item\emph{Monotonicity}: i.e. $x\leq y$ implies $\rho(y)\leq\rho(x)$;
\item\emph{Cash-invariance}: i.e. $\rho(x+r)=\rho(x)-r$ for all $r\in\R$.
\end{enumerate}
Such a function $\rho$ is called a \emph{convex risk measure}. 
The notion of convex risk measure was independently introduced in \cite{follmer2002convex} and \cite{fritelli03} as a generalization of the notion of \emph{coherent risk measure} introduced in \cite{artzner1999coherent}.

Associated to the model space $\mathcal{X}$, we can consider a dual pair. 
Namely, the K\"{o}the dual space of $\mathcal{X}$ is defined by
\[
\mathcal{X}^\#:=\{y\in L^0(\EE)\colon x y\in L^1(\EE)\mbox{ for all }x\in\mathcal{X}\}.
\] 
Then $\mathcal{X}^\#$ is also a solid subspace of $L^1(\EE)$ with $\R\subset\mathcal{X}^\#$.  
This gives rise to the dual pair $\langle\mathcal{X},\mathcal{X}^\#\rangle$ associated to the bilinear form $(x,y)\mapsto\E_\PP[x y]$ and the weak topologies $\sigma(\mathcal{X},\mathcal{X}^\#)$ and $\sigma(\mathcal{X}^\#,\mathcal{X})$.\\

If $\mathcal{X}=L^p(\mathcal{E})$ with $1\leq p\leq \infty$, it is known that $\mathcal{X}^\#=L^q(\mathcal{E})$, where $q:=(1-1/p)^{-1}$ is the H\"{o}lder conjugate of $p$, see e.g.~\cite[Example 29.4]{zaanen2012introduction}. 
Suppose that $\phi$ is a Young function and let $\psi(r):=\sup_{s\geq 0}\{r s - \phi(s)\}$ be the conjugate Young function of $\phi$.  
If $\mathcal{X}=L^\phi(\mathcal{E})$, then one has that $\mathcal{X}^\#=L^\psi(\mathcal{E})$, see e.g.~\cite{luxemburg44701banach, zaanen1983riesz}.  
If  $\mathcal{X}=H^\phi(\mathcal{E})$ and $\phi$ is finite-valued (otherwise $H^\phi(\mathcal{E})=\{0\}$), then $\mathcal{X}^\#=L^\psi(\mathcal{E})$, see e.g.~\cite{owari2014lebesgue}.\\

The \emph{Fenchel transform} of a convex risk measure $\rho$ is defined to be
\[
\rho^\#(y):=\sup\{\E_\PP[x y]-\rho(x) \colon x\in\mathcal{X}\}.
\]

Duality theory of convex risk measures is aimed to study when the Fenchel transform is involutive.  
More precisely, given a convex risk measure $\rho$ say that:
\begin{itemize}
\item $\rho$ is \emph{representable} if it admits the following dual representation:
\[
\rho(x)=\sup\{\E_\PP[x y]-\rho^\#(y)\colon y\in\mathcal{X}^\#\}\quad\text{ for all }x\in\mathcal{X}.
\]
\item $\rho$ attains its representation whenever for any $x\in\mathcal{X}$ there exists a $y\in \mathcal{X}^\#$ such that
\[
\rho(x)=\E_\PP[x y]-\rho^\#(y).
\]
\end{itemize} 

\begin{rem}
We have that $\rho^\#(y)<\infty$ only if $y\leq 0$ and $\E[y]=-1$.\footnote{Indeed, suppose that $y\in X^\#$ and fix $n\in\N$. 
Then $\rho^\#(y)\geq \E[n 1_{\{y\geq 0\}}y] - \rho(n 1_{\{y\geq 0\}})\geq n\E[y^+] -\rho(0).$ 
Since $n$ is arbitrary, $\rho^\#(y)<\infty$ only if $y^+=0$. 
We also have that $\rho^\#(y)\geq \E[n y] - \rho(n)\geq n(\E[y]+1) -\rho(0)$. 
Being $n$ arbitrary, we conclude that $\rho^\#(y)<\infty$ only if $\E[y]=-1$.}
 Thus $\rho$ is representable if and only if 
\begin{equation}
\label{eq:stressTest}
\rho(x)=\sup\{\E_\PP[x y]-\rho^\#(y)\colon y\in\mathcal{X}^\#,\:y\leq 0,\E_\PP[y]=-1\}\quad\text{ for all }x\in\mathcal{X}.
\end{equation}
Notice that an element $y\in L^1(\EE)$ with $y\leq 0$ and $\E_\PP[y]=-1$ can be identified with a probability measure $Q_y\ll\PP$ via the Radon-Nikodym derivative $y=-\frac{d Q_y}{d \PP}$. 
The economic interpretation of the representation (\ref{eq:stressTest}) is that a convex risk measure can be seen as a stress test of the financial position $x$ among the different market models given by the probabilities $Q_y$ and the penalty function $\rho^\#$. 
\end{rem}

Next, we recall some properties that matter in duality theory of convex risk measures:

Let $f$ be a function from $\mathcal{X}$ to the extended real numbers $\overline{\mathbb{R}}$. 
For any $r\in\mathbb{R}$, we define the sublevel set 
\[
V_r(f):=\left\{x\in\mathcal{X}\colon f(x)\leq r\right\}.
\]
Say that: 
\begin{itemize}
\item $f$ has the \emph{Fatou property} if 
\[
\lim_n x_n=x\text{ a.s., }y\in\mathcal{X},\:|x_n|\leq y\text{ for all }n\in\N\text{ implies }\liminf_n f(x_n)\geq f(x);
\]
\item $f$ has the \emph{Lebesgue property} if 
\[
\lim_n x_n=x\text{ a.s., }y\in\mathcal{X},\:|x_n|\leq y\text{ for all }n\in\N\text{ implies }\lim_n f(x_n)=f(x);
\]
\item  $f$ is  \emph{law invariant} if $f(x)=f(y)$ whenever $x$ and $y$ have the same law (i.e.~$\PP(x\leq r)=\PP(y\leq r)$ for each $r\in\R$);
\item $f$ is \emph{lower semicontinuous} w.r.t.~$\sigma(\mathcal{X},\mathcal{X}^\#)$, 
if $V_r(f)$ is closed  w.r.t. $\sigma(\mathcal{X},\mathcal{X}^\#)$ for each $r\in\R$;
\item $f$ is  \emph{inf-compact} w.r.t. $\sigma(\mathcal{X},\mathcal{X}^\#)$, if $V_r(f)$ is compact  w.r.t. $\sigma(\mathcal{X},\mathcal{X}^\#)$ for each $r\in\R$.
\end{itemize}

\subsection{Multi-period setup: conditional risk measures}

The notion of conditional risk measure was independently introduced by \cite{detlefsen2005conditional} and \cite{bionnadal04}. 
Next, we recall the main elements of duality theory of conditional risk measures. 
Namely, we adopt the module-based approach introduced in~\cite[Section 3]{kupper11}.

Now, we suppose that $\F$ is a sub-$\sigma$-algebra of $\mathcal{E}$, which models the available market information at some future date $t\in (0,T)$. 
Let us introduce some notation. 
We denote by $L^0_{+}(\F)$, $L^0_{++}(\F)$, and $\overline{L^0}(\F)$ the spaces of (classes of equivalence of) $\F$-measurable random variables with values in the intervals $[0,\infty)$, $(0,\infty)$, and $[-\infty,\infty]$, respectively. 

Let $\bar\F$ be the probability algebra associated to $\left(\Omega,\mathcal{F},\PP\right)$, where $\bar\F$ is defined by identifying events modulo null sets.  
It is well-known that $\bar\F$ is a complete Boolean algebra which satisfies the countable chain condition (ccc), i.e. every family of positive pairwise disjoint elements in $\bar\F$ is at most countable. 
The $0$ of $\bar\F$ is represented by the empty set $\emptyset$ and the unity $I$ of $\bar\F$ is represented by $\Omega$. 
We denote by $p(I)$ the set of all partitions of $I$ to $\bar\F$.  
Given $a\in\bar\F$, we write $1_a$  for the class in $L^0(\F)$ of the characteristic function $1_A$ of some representative $A\in\F$ of $a$.
Given a partition $(a_k)_{k\in\mathbb{N}}\in p(I)$ and a sequence $(x_k)_{k\in\mathbb{N}}$, we define $\sum 1_{a_k}x_k:=\lim_k 	 \sum_{i=1}^k 1_{a_i}x_i$ a.s.. 
  
Classically, the conditional expectation $\mathbb{E}[\cdot|\F]$ is defined for elements with finite expectation. 
We consider the extended conditional expectation. 
Namely, suppose that $x\in L^0(\mathcal{E})$ satisfies that 
 at least one of $\lim_n \E_\PP[x^+\wedge n|\F]$ and $\lim_n \E_\PP[x^-\wedge n|\F]$ (a.s.) is finite,  
then we define the  \emph{extended conditional expectation}  of $x$ to be
\[
\E_\PP[x|\F]:=\lim_n \E_\PP[x^+\wedge n|\F]-\lim_n \E_\PP[x^-\wedge n|\F]\in \overline{L^0}(\F). 
\]

Now, our model space is an $L^0(\F)$-submodule $\mathscr{X}$ of $L^0(\mathcal{E})$ which satisfies the following properties:
\begin{itemize}
\item $\mathscr{X}$ is solid;
\item $\E_\PP[|x||\F]<\infty$ for all $x\in\mathscr{X}$;
\item $L^0(\F)\subset \mathscr{X}$;
\item $\mathscr{X}$ is \emph{stable}, that is, $\sum 1_{a_k}x_k\in \mathscr{X}$ whenever $(a_k)\in p(I)$ and $(x_k)\subset \mathscr{X}$.
\end{itemize}

\begin{example}

The following $L^0(\F)$-submodules of $L^0(\mathcal{E})$ satisfy the properties above:

\begin{enumerate}
\item \emph{$L^p$ type modules} (see \cite{kupper03}): 
We define
\[
L^\infty_\F(\mathcal{E}):=\left\{  x\in L^0(\mathcal{E}) \colon  |x|\leq\eta,\:\text{for some }\eta\in L^0(\F)\right\},
\]

if $1\leq p<\infty$, let
\[
L^p_\F(\mathcal{E}):=\left\{ x\in L^0(\mathcal{E}) \colon \E_\PP\left[\left|x\right|^{p}|\mathcal{F}\right]<\infty \right\}.
\]
\item \emph{Orlicz type modules} (see \cite{vogelpoth20090}) and \emph{Orlicz-heart type modules}: Let $\phi:[0,\infty)\rightarrow[0,\infty]$ be a Young function  and let
\[
L^\phi_\F(\mathcal{E}):=\{x\in L^0(\mathcal{E}) \colon \exists \eta\in L^0_{++}(\F),\:\E_\PP[\phi(\eta|x|)|\F]\in L^0(\F) \}\quad\mbox{ (Orlicz type module),}
\]
\[
H^\phi_\F(\mathcal{E}):=\{x\in L^0(\mathcal{E}) \colon \forall \eta\in L^0_{++}(\F),\: \E_\PP[\phi(\eta|x|)|\F]\in L^0(\F) \}\quad\mbox{ (Orlicz-heart type module)}.\\
\]

\end{enumerate}

\end{example}

Our model module $\mathscr{X}$ is going to describe all possible final payoffs of the positions at $T$.
 
The riskiness at time $t$ of any financial position $x\in\mathscr{X}$ is uncertain and contingent to the information encoded in $\F$. 
Thus the riskiness is quantified by a function $\rho:\mathscr{X}\rightarrow L^0(\F)$ which satisfies:
\begin{enumerate}
\item \emph{$L^0(\F)$-convexity}: $\rho(\eta x + (1-\eta)y)\leq \eta\rho(x) + (1-\eta)\rho(y)$ whenever $\eta\in L^0(\F)$ with $0\leq \eta\leq 1$ and $x,y\in\mathscr{X}$;
\item \emph{Monotonicity}: if $x\leq y$ in $\mathscr{X}$, then $\rho(y)\leq\rho(x)$;
\item \emph{$L^0(\F)$-cash invariance}: $\rho(x+\eta)=\rho(x)-\eta$ whenever $\eta\in L^0(\F)$, $x\in\mathscr{X}$.
\end{enumerate}
Such a function is called a \emph{conditional risk measure}.

Dual systems of modules were introduced and studied in~\cite{kusraev1985vector}. 
Associated to the model space $\mathscr{X}$, we can consider a dual system of $L^0(\F)$-modules. 
Namely, we define the \emph{Köthe dual $L^0(\F)$-module of $\mathscr{X}$} to be
\[
\mathscr{X}^\#:=\left\{ y\in L^0(\mathcal{E}) \colon x y\in L^1_\F(\mathcal{E})\text{ for all }x\in\mathscr{X}\right\}.
\]

It is simple to verify that $\mathscr{X}^\#$ enjoys the same properties as  $\mathscr{X}$;  
namely, $\mathscr{X}^\#$ is a solid and stable $L^0(\F)$-submodule with $L^0(\F)\subset\mathscr{X}^\#\subset L^1_\F(\mathcal{E})$.  

The dual system $\langle\mathscr{X},\mathscr{X}^\#\rangle$ allows for the definition of the following module analogue of the Fenchel transform:
\[
\rho^\#(y):=\esssup\{\E_\PP[x y|\F]-\rho(x) \colon x\in \mathscr{X}\} \quad \mbox{ for }y\in\mathscr{X}^\#. 
\] 

Again, we are interested in the involutivity of the Fenchel transform. 
Thus we introduce the following nomenclature: 
Given a conditional risk measure  $\rho:\mathscr{X}\rightarrow L^0(\F)$, say that:
\begin{itemize}
\item $\rho$ is \emph{representable} if  
\[
\rho(x)=\esssup\{\E_\PP[x y|\F]-\rho^\#(y)\colon y\in\mathscr{X}^\#\}\quad\mbox{ for all }x\in\mathscr{X}.
\]
\item $\rho$ attains its representation if for any $x\in\mathscr{X}$ there exists $y\in\mathscr{X}^\#$ such that
\[
\rho(x)=\E_\PP[x y|\F]-\rho^\#(y).
\]
\end{itemize} 

\begin{rem}
Due to \cite[Corollary 3.14]{kupper11}, a conditional risk measure $\rho$ is representable if and only if
\[
\rho(x)=\esssup\{\E_\PP[x y|\F]-\rho^\#(y)\colon y\in\mathscr{X}^\#,\:y\leq 0,\: \mathbb{E}_{\PP}[y|\F]=-1\}\quad\mbox{ for all }x\in\mathscr{X}.
\]
\end{rem}

Next, we will introduce some notions that are useful in the dual representation of a conditional risk measure. 

Given the dual system of $L^0(\F)$-modules $\langle\mathscr{X},\mathscr{X}^\#\rangle$, we can define the so-called \emph{stable weak topologies} induced by $\langle\mathscr{X},\mathscr{X}^\#\rangle$. 
Namely, given a partition $(a_k)\in p(I)$, a family $(F_k)$ of non-empty finite subsets of $\mathscr{X}^\#$, and $\varepsilon\in L^0_{++}(\F)$, we define
\[
U_{(F_k),(a_k),\varepsilon}:=\left\{ x\in\mathscr{X} \colon \sum 1_{a_k}\underset{y\in F_k}\esssup|\E_\PP[x y|\F]|<\varepsilon \right\}.
\]
The collection of sets $$\mathscr{B}_{\sigma_s(\mathscr{X},\mathscr{X}^\#)}:=\left\{ x+U_{(F_k),(a_k),\varepsilon}\colon x\in\mathscr{X},\:(a_k)\in p(I),\:\emptyset\neq F_k\subset\mathscr{X}^\#\text{ finite, }\varepsilon\in L^0_{++}(\F)\right\}$$
is a base for a topology on $\mathscr{X}$, which  will be denoted by $\sigma_s(\mathscr{X},\mathscr{X}^\#)$. 
Similarly, we define $\sigma_s(\mathscr{X}^\#,\mathscr{X})$. 

Stable weak topologies were introduced in~\cite[1.1.8]{kusraev1985vector} as the topology induced by the multi-norm associated to a dual system of modules. 
Also, they are the transcription in a modular setting of the conditional weak topologies introduced in~\cite{DJKK13}.

Suppose that   $f$ is a function from $\mathscr{X}$ to $\overline{L^0}(\F)$. 
For any $\eta\in L^0(\F)$, we define the sublevel set
\[
V_\eta(f):=\left\{x\in\mathscr{X}\colon f(x)\leq\eta\right\}.
\]
Say that:
\begin{itemize}
\item $f$ has the \emph{Fatou property} if 
\[
\lim_n x_n=x\text{ a.s., }y\in\mathscr{X},\:|x_n|\leq y\text{ for all }n\in\N\text{ implies }\liminf_n f(x_n)\geq f(x);
\]
\item $f$ has the \emph{Lebesgue property} if 
\[
\lim_n x_n=x\text{ a.s., }y\in\mathscr{X},\:|x_n|\leq y\text{ for all }n\in\N\text{ implies }\lim_n f(x_n)=f(x)\text{ a.s.};
\]	
	
\item $f$ is \emph{conditionally law invariant} is $f(x)=f(y)$ whenever $x$ and $y$ have the same \emph{conditional law} (i.e. $\PP(x\leq\eta|\F)=\PP(y\leq\eta|\F)$ for each $\eta\in L^0(\F)$);

\item lower semicontinuous w.r.t. $\sigma_s(\mathscr{X},\mathscr{X}^\#)$ if $V_\eta(f)$ is closed w.r.t. $\sigma_s(\mathscr{X},\mathscr{X}^\#)$ for every $\eta\in L^0(\F)$.	
\end{itemize} 

Next, we will recall some notions that we will be needed later. 
\begin{itemize}
\item A non-empty subset $S$ of $L^0(\mathcal{E})$ is \emph{stable} if $\sum 1_{a_k}x_k\in S$ whenever $(x_k)\subset S$ and $(a_k)\in p(I)$;
\item A collection $\mathscr{B}$ of non-empty subsets of $L^0(\mathcal{E})$ is said to be \emph{stable} if for any sequence $(S_k)$ of members of $\mathscr{B}$ and any partition $(a_k)\in p(I)$ one has
\[
\sum 1_{a_k}S_k=\left\{ \sum_k 1_{a_k}x_k \colon  x_k\in S_k\text{ for all }k\right\}\in \mathscr{B};
\]
\item A \emph{stable filter base} is a filter base $\mathscr{B}$ on $L^0(\mathcal{E})$ such that $\mathscr{B}$ is a stable collection consisting of stable subsets of $L^0(\mathcal{E})$; 
\item
A non-empty subset $S$ of $\mathscr{X}$ is \emph{stably compact} with respect to $\sigma_s(\mathscr{X},\mathscr{X}^\#)$, if $S$ is stable and any stable filter base $\mathscr{B}$ on $S$ has a cluster point $x\in S$ w.r.t.~$\sigma_s(\mathscr{X},\mathscr{X}^\#)$;
\item  A function $f:\mathscr{X}\to \overline{L^0}(\mathcal{E})$ is said to be \emph{stably inf-compact} w.r.t. $\sigma_s(\mathscr{X},\mathscr{X}^\#)$ if $V_\eta(f)$ is stably compact w.r.t.~$\sigma_s(\mathscr{X},\mathscr{X}^\#)$ for every $\eta\in L^0(\F)$ such that $V_\eta(f)\neq\emptyset$.
\end{itemize}

The notion of stability is crucial in some related frameworks. 
In Boolean valued analysis it is used the terminology \emph{cyclic} or \emph{universally complete} $\A$-sets (here $\A$ is any complete Boolean algebra, for instance we can take $\A=\bar{\F}$), see \cite{kusraev2012boolean}. 
In particular, in the case of dual systems of modules this notion was introduced in~\cite{kusraev1985vector}.
 
In conditional set theory it is used the terminology stable set and stable collection, see \cite{DJKK13}. 
Actually, the notion of conditional set is a reformulation of that of cyclic $\A$-set. 
However, it should be mentioned that conditional set theory provides us with an intuitive and useful tool for dealing with $\A$-sets and their Boolean valued representation. 
In theory of $L^0$-modules the notion of stability is called the \emph{countable concatenation property}, see~\cite{guo10}. 

Stable compactness was first time studied by Kusraev \cite{cyclicCompactness} under the name of \emph{cyclic compactness}. 
The notion of stable filter base and stable compactness were defined in \cite{DJKK13}. 
The transcriptions of these notions in $L^0$-modules is studied in \cite{L0compactness}.

\section{Foundations of Boolean valued models}

The precise formulation of Boolean valued models requires some familiarity with the basics of set theory and logic, and in particular with first-order logic, ordinals and transfinite induction. 
For the convenience of the reader, we will give some background of this theory. 
For a more detailed description we refer the reader to \cite{kusraev2012boolean}.

Let us consider a universe of sets $V$ satisfying the axioms of the Zermelo-Fraenkel set theory with the axiom of choice (ZFC), and  a first-order language $\mathcal{L}$, which allows for the formulation of statements about the elements of $V$. 
In the universe $V$ we have all possible mathematical objects (real numbers, topological spaces, and so on). 
The language $\mathcal{L}$ consists of names for the elements of $V$ together with a finite list of symbols for logic symbols ($\forall$, $\wedge$, $\neg$ and parenthesis), variables and the predicates  $=$ and $\in$.  Though we usually use a much richer language by introducing more and more intricate definitions, in the end any usual mathematical statement can be written using only those mentioned. The elements of the universe $V$ are classified into a transfinite hierarchy: $V_0\subset V_1\subset V_2 \subset \cdots V_\omega \subset V_{\omega+1}\subset \cdots$, where $V_0 = \emptyset$, $V_{\alpha+1} = \mathcal{P}(V_\alpha)$ is the family of all sets whose elements come from $V_\alpha$, and $V_\beta = \bigcup_{\alpha<\beta}V_\alpha$ for limit ordinal $\beta$.

The following constructions and principles work for any complete Boolean algebra $\mathcal{A}$, even if it is not associated to a probability space or even does not have the countable chain condition. 
However, for the sake of simplicity, we will consider our underlying probability algebra $\mathcal{A}:=\bar{\F}$, which encodes the 
future market information.

We will construct $\VA$, the \emph{Boolean valued model} of $\A$, whose elements we interpret as objects which we can talk about at the future time $t$. 
 We proceed by induction over the class $\Ord$ of ordinals of the universe $V$. 
 We start by defining $\VA_0:=\emptyset$.  
 If $\alpha+1$ is the successor of the ordinal $\alpha$, we define
 $$\VA_{\alpha+1}:=\left\{ u \colon u\text{ is an $\A$-valued function with }\dom(u)\subset\VA_\alpha  \right\}.$$
 If  $\alpha$ is a limit ordinal  $V_\alpha^{(\A)}:=\underset{\xi<\alpha}\bigcup V_{\xi}^{(\A)}$. 
 Finally, let $V^{(\A)}:=\underset{\alpha\in Ord}\bigcup V_{\alpha}^{(\A)}.$
 
 The idea is that any member $v$ of the class  $V^{(\A)}$ is a \emph{fuzzy} set in the sense that, for $v\in \text{dom}(u)$, $v$ will become an element of $u$ at the future time $t$ if $u(v)$ happens. 
 Given  $u$ in $V^{(\A)}$, we define its \emph{rank} as the least ordinal $\alpha$ such that $u$ is in $\VA_{\alpha+1}$.
 
 We consider a first-order language which allows us to produce statements about $\VA$.
 Namely, let $\mathcal{L}^{(\A)}$ be the first-order language which is the extension of $\mathcal{L}$ by adding \emph{names} for each element in $\VA$. 
 Throughout, we will not distinguish between an element in $V^{(\A)}$ and its name in $\mathcal{L}^{(\A)}$. 
 Thus, hereafter, the members of $\VA$ will be referred to as names.  
  
 Suppose that $\varphi$ is a formula in set theory, that is, $\varphi$ is constructed  by applying logical symbols to
atomic formulas $u=v$ and $u\in v$.  
 If $\varphi$ does not have any free variable and all the constants in $\varphi$ are names in $V^{(\A)}$, then we define its \emph{Boolean truth value}, say $\llbracket\varphi\rrbracket$, which is a member of $\A$ and is constructed by induction in the length of $\varphi$ by naturally giving Boolean meaning to the predicates  $=$ and $\in$, the logical connectives and the quantifiers.
 
We start by defining the Boolean truth value of the \emph{atomic formulas} $u\in v$ and $u=v$ for $u$ and $v$ in $\VA$. Namely, proceeding by transfinite recursion we define
 \[
 \llbracket u\in v\rrbracket=\underset{t\in\dom(v)} \bigvee v(t)\wedge\llbracket t=u\rrbracket,
 \]
 \[
 \llbracket u=v\rrbracket=\underset{t\in\dom(u)}\bigwedge \left(u(t)\Rightarrow \llbracket t\in v\rrbracket\right) \wedge \underset{t\in\dom(v)}\bigwedge \left(v(t)\Rightarrow \llbracket t\in u\rrbracket\right),
 \]

where, for $a,b\in\A$, we denote $a\Rightarrow b:=a^c\vee b$.
 For non-atomic formulas we have 
 \[
 \llbracket (\exists x)\varphi(x)\rrbracket:=\underset{u\in \VA}\bigvee \llbracket \varphi(u)\rrbracket\quad\mbox{ and }\quad \llbracket (\forall x)\varphi(x)\rrbracket:=\underset{u\in \VA}\bigwedge \llbracket \varphi(u)\rrbracket;
 \]
\[
\llbracket \varphi \vee \psi\rrbracket:=\llbracket \varphi\rrbracket \vee \llbracket\psi\rrbracket, \quad\llbracket \varphi \wedge \psi\rrbracket:=\llbracket \varphi\rrbracket \wedge \llbracket\psi\rrbracket,\quad\llbracket \varphi \Rightarrow \psi\rrbracket:=\llbracket \varphi\rrbracket^c \vee \llbracket\psi\rrbracket \quad\mbox{ and }\quad \llbracket \neg\varphi\rrbracket:=\llbracket\varphi\rrbracket^c.
\]
  
We say that a formula $\varphi$ is satisfied within $V^{(\A)}$, and write $V^{(\A)}\models \varphi$, whenever it is true with the Boolean truth value, that is, $\llbracket\varphi \rrbracket=I$.

We say that two names $u,v$ are equivalent when $\llbracket u=v \rrbracket =I$. 
It is not difficult to verify that the Boolean truth value of a formula is not affected when we change a name by an equivalent one. 
However, the relation $\llbracket u=v \rrbracket =I$ does not mean that the functions $u$ and $v$ (considered as elements of $V$) coincide. 
For example, the empty function $u:=\emptyset$ and the function $v:\{\emptyset\}\rightarrow \A$ $v(\emptyset):=0$ are different as functions; however, $\llbracket u=v\rrbracket=I$.  
In order to avoid technical difficulties, we will consider the so-called \emph{separated universe}. 
Namely, let $\VAA$ be the subclass of $\VA$ defined by choosing a representative of the least rank in each class of the equivalence relation $\{(u,v)  \colon \llbracket u=v\rrbracket=I \}$.

The universe $V$ can be embedded into $\VA$. 
Given a set $x$ in $V$, we define its canonical name $\check{x}$ in $\VA$ by transfinite induction. 
Namely, we put $\check{\emptyset}:=\emptyset$ and for $x$ in $V$ we define $\check{x}$ to be the unique representative in $\VAA$ of the name given by the function
\[
\left\{\check{y}\colon y\in x \right\}\rightarrow \A,\quad \check{y}\mapsto I.
\]

Given a name $u$ with $\llbracket u\neq\emptyset\rrbracket=I$ we define its \emph{descent} by  
\[
u{\downarrow} := \{ v\in\VAA \colon \llbracket v\in u\rrbracket=I\}.
\]
$V^{(\mathcal{A})}$ is a model of ZFC. More precisely we have:

\begin{thm}(Transfer Principle)
\label{thm: Transfer}
If $\varphi$ is a theorem of ZFC, then $\VA\models \varphi$.
\end{thm}

Other two important principles are the following:

\begin{thm}(Maximum Principle)
\label{thm: Maximum}
Let $\varphi(x_1,\ldots,x_n)$ be a formula with free variables $x_1,\ldots,x_n$. 
Then there exist names $u_1,\ldots,u_n$ such that $\llbracket\varphi(u_1,\ldots,u_n)\rrbracket=\llbracket (\exists x_1) \ldots (\exists x_n)\varphi(x_1,\ldots,x_n)\rrbracket$. 

\end{thm}

\begin{thm}(Mixing Principle)
\label{thm: Mixing}
Let $(a_k)\in p(I)$ and let $(u_k)$ be a sequence of names. 
Then there exists a unique member $u$ of $\VAA$ such that $\llbracket u=u_k\rrbracket\geq a_k$ for all $k\in\N$. 
\end{thm}

Given a partition $(a_k)\in p(I)$ and  a sequence $(u_k)$ of elements of $\VA$, we denote by $\sum u_k a_k$, the unique name $u$ in $\VAA$ satisfying $\llbracket u=u_k\rrbracket\geq a_k$ for all $k\in\N$. \\

The following result is very useful to manipulate Boolean truth values:

\begin{prop}\label{thm: forallExists}
Let $\varphi(x)$ be a formula with a free variable $x$ and  $v$ a name with $\llbracket v\neq\emptyset\rrbracket=I$. 
Then:
\[
\llbracket (\forall x \in v)\varphi(x)\rrbracket=\underset{u\in v{\downarrow}}\bigwedge\llbracket\varphi(u) \rrbracket,\quad
\llbracket (\exists x \in v)\varphi(x)\rrbracket=\underset{u\in v{\downarrow}}\bigvee\llbracket\varphi(u) \rrbracket.
\]
Moreover, one has
\begin{enumerate}
\item $\llbracket(\forall x \in v)\varphi(x)\rrbracket=I$ if and only if $\llbracket\varphi(u)\rrbracket=I$ for all $u\in v{\downarrow}$;
\item $\llbracket(\exists x \in v)\varphi(x)\rrbracket=I$ if and only if there exists $u\in{v}{\downarrow}$ such that $\llbracket\varphi(u)\rrbracket=I$.
\end{enumerate}
\end{prop}

In general, in the universe $V^{(\mathcal{A})}$ we have all possible mathematical objects (real numbers, topological spaces, and so on). 
If $u$ is a name which satisfies $\llbracket u\textnormal{ is a function}\rrbracket=I$, that is, $u$ satisfies the definition of function in the language $\mathcal{L}^{(\A)}$, then we say that $u$ is a \emph{name for a function}. 
Of course, this can be done for any mathematical concept. Thus, in the sequel of this article, we will use the terminology \emph{name for a vector space}, \emph{name for a topology}, and so on without further explanations.

\begin{defn}
Suppose that $u,v$ are two names with $\llbracket (u\neq \emptyset)\wedge (v\neq\emptyset)\rrbracket=I$.  
A function $f:u{\downarrow}\rightarrow v{\downarrow}$ such that 
$$\llbracket w=t\rrbracket\leq \llbracket f(w)=f(t)\rrbracket\quad\mbox{ for all }w,t\in u{\downarrow}$$
is called \emph{extensional}. 
\end{defn}

Extensional functions allows for the definition of names for functions. 
More precisely, we have the following:

\begin{prop}
\label{thm: extensional}
Let $u,v$ be names with $\llbracket (u\neq \emptyset)\wedge (v\neq\emptyset)\rrbracket=I$ and suppose that $f:u{\downarrow}\rightarrow v{\downarrow}$ is an extensional function.  
Then there exists a name $f{\auu}$ for a function from $u$ to $v$, 
such that $\llbracket f{\auu}(t)=f(t)\rrbracket=I$ for all $t\in u{\downarrow}$.  
\end{prop}

\section{A transfer principle between duality theory of convex risk measure and duality theory of conditional risk measures}

Let us go back to our model probability space $(\Omega,\mathcal{E},\PP)$ with $\F\subset\mathcal{E}$. 
Next, we state the main result of the present paper, which allows for the interpretation of a conditional risk measure $\rho:\mathscr{X}\rightarrow L^0(\F)$ as a name for a convex risk measure, let us say $\rho{\uparrow}$, defined on some space of random variables, and relates the properties of $\rho$ with the properties of $\rho{\uparrow}$ in the set-theoretic model $\VA$. 
In other words, this result establishes a transfer principle between duality theory of convex risk measures and duality theory of conditional risk measures.

\begin{thm}\label{prop: transferI}
Let $\rho:\mathscr{X}\to L^0(\mathcal{F})$ be a conditional risk measure. 
Then there exist members $\rho\uparrow$  and $\mathscr{X}{\uparrow}$ of $\VA$ such that
\begin{align*}
V^{(\A)}\models&\text{ there exists a probability space $(X,\Sigma,Q)$ such that,}\\
               &\mathscr{X}{\uparrow}\text{ is a solid subspace of $L^1(\Sigma)$ with $\mathbb{R}\subset \mathscr{X}{\uparrow}$},\\
               &\text{and $\rho{\uparrow}:\mathscr{X}{\uparrow}\to\mathbb{R}$  is a convex risk measure},
\end{align*}
 and so that the names $\rho\uparrow$  and $\mathscr{X}{\uparrow}$ satisfy the following:
\begin{enumerate}

\item $\rho$ is representable iff $\llbracket\rho{\uparrow}\textnormal{ is representable}\rrbracket=I$.
\item $\rho$ attains its representation iff 
$\llbracket\rho{\uparrow}\textnormal{ attains its representation}\rrbracket=I$.
\item $\rho$ has the Fatou property iff $\llbracket\rho{\uparrow}\textnormal{ has the Fatou property}\rrbracket=I$.
\item $\rho$ has the Lebesgue property iff $\llbracket\rho{\uparrow}\textnormal{ has the Lebesgue property}\rrbracket=I$.
\item $\rho$ is conditionally law invariant iff $\llbracket\rho{\uparrow}\textnormal{ is law invariant}\rrbracket=I$.

\item $\rho$ is lower semicontinuous  w.r.t $\sigma_s(\mathscr{X},\mathscr{X}^\#)$ iff 

$\llbracket\rho{\uparrow}\textnormal{ is lower semicontinuous w.r.t. }\sigma(\mathscr{X}{\uparrow},\mathscr{X}{\uparrow}^\#)\rrbracket=I$.
\item $\rho^\#$ is stably inf-compact w.r.t $\sigma_s(\mathscr{X}^\#,\mathscr{X})$ iff 

$\llbracket\rho{\uparrow}^\#\textnormal{ is inf-compact w.r.r. }\sigma(\mathscr{X}{\uparrow}^\#,\mathscr{X}{\uparrow})\rrbracket=I$.
\item If $\mathscr{X}=L^p_\F(\mathcal{E})$ with $1\leq p\leq \infty$, then  $\llbracket \mathscr{X}{\uparrow}=L^{p}(\Sigma)\rrbracket=I$. 
In that case,  $\mathscr{X}^\#=L^q_\F(\mathcal{E})$ where $q$ is the H\"{o}lder conjugate of $p$.
\item If $\mathscr{X}=L^\phi_\F(\mathcal{E})$ with $\phi$ a Young function, 
then there is a name $\tilde{\phi}$ for a Young function such that 
$\llbracket \mathscr{X}{\uparrow}=L^{\tilde{\phi}}(\Sigma)\rrbracket=I$. 
In that case, $\mathscr{X}^\#=L^\psi_\F(\mathcal{E})$ where $\psi$ is the conjugate Young function of $\phi$.
\item If $\mathscr{X}=H^\phi_\F(\mathcal{E})$ with $\phi$ a finite-valued Young function, 
then there is a name $\tilde{\phi}$ for a finite-valued Young function such that 
$\llbracket \mathscr{X}{\uparrow}=H^{\tilde{\phi}}(\Sigma)\rrbracket=I$. 
In that case, $\mathscr{X}^\#=L^\psi_\F(\mathcal{E})$ where $\psi$ is the conjugate Young function of $\phi$.
\end{enumerate}
\end{thm}

The  proof of the theorem above is postponed to next section. 
Instead, we focus first on some instances of application.

Theorem \ref{prop: transferI} together with the transfer principle of Boolean-valued models allow for the interpretation of well-known results about the dual representation of convex risk measures as new theorems about the dual representation of conditional risk measures. 

For example, suppose that $\rho:\mathcal{X}\rightarrow\R$ is a convex risk measure. 
As a consequence of the Fenchel-Moreau theorem (see \cite[Theorem 2.1]{kaina2009convex}) applied to the weak topology $\sigma(X,X^\#)$ we have that $\rho$ is representable if and only if $\rho$ is lower semicontinuous w.r.t. $\sigma(X,X^\#)$. 

Moreover, we have the following dual representation result:

\begin{thm}\label{thm: owari}\cite[Theorem 1.1]{owari2014lebesgue}
Let $\rho:\mathcal{X}\rightarrow\R$ be a convex risk measure. 
Then $\rho$ is representable 
if and only if $\rho$ is lower semicontinuous w.r.t. $\sigma(\mathcal{X},\mathcal{X}^\#)$. 
In that case, the following statements are equivalent:
\begin{enumerate}
\item $\rho$ attains its representation;
\item $\rho$ has the Lebesgue property;
\item $\rho^\#$ is inf-compact w.r.t. $\sigma(\mathcal{X}^\#,\mathcal{X})$.
\end{enumerate}
\end{thm}

Let $\varphi$ denote the  theorem above. 
Due to the transfer principle, it is satisfied that $\llbracket\varphi\rrbracket=I$. 
In view of Theorem \ref{prop: transferI}, we have that the statement below is just a reformulation of $\llbracket\varphi\rrbracket=I$, so no proof is needed.

\begin{thm}\label{thm: rep}
Let $\rho:\mathscr{X}\rightarrow L^0(\F)$ be a conditional risk measure. 
Then $\rho$ is representable, i.e. 
\[
\rho(x)=\esssup\{\E_\PP[x y|\F]-\rho^\#(y)\colon y\in\mathscr{X}^\#,\:y\leq 0,\E_\PP[y|\F]=-1\}\quad\text{ for all }x\in\mathscr{X}
\]
if and only if $\rho$ is lower semicontinuous w.r.t. $\sigma_s(\mathscr{X},\mathscr{X}^\#)$. 

In that case, the following are equivalent:

\begin{enumerate}
\item $\rho$ attains its representation, i.e. for every $x\in \mathscr{X}$ there exists $y\in \mathscr{X}^\#$ with $y\leq 0$ and $\E[y|\F]=-1$ such that $
	\rho(x)=\E_\PP[x y|\F]-\rho^\#(y);$
	 
	\item $\rho$ has the Lebesgue property;
	\item $\rho^\#$ is stably inf-compact w.r.t. $\sigma_s(\mathscr{X}^\#,\mathscr{X})$.
\end{enumerate}
\end{thm}

Suppose that $\mathcal{X}=L^\infty(\mathcal{E})$. 
Then, the so-called Jouini-Schachermayer-Touzi theorem (see~\cite[Theorem 2]{delbaen2009differentiability} and for its original form see~\cite{JST06}) asserts that in Theorem \ref{thm: owari} we can replace the lower semicontinuity by the Fatou property.\footnote{Actually, the Fatou property is automatically satisfied when $\rho$ is law invariant, see  \cite{JST06}.}
Thus, the transfer principle together with Theorem \ref{prop: transferI} yields the following:

\begin{thm}
\label{thm: repII}
Let $\rho:L^\infty_\F(\mathcal{E})\rightarrow L^0(\F)$ be a conditional risk measure. 
Then  $\rho$ has the Fatou property if and only if it admits a representation
\[
\rho(x)=\esssup\{\E_\PP[x y|\F]-\rho^\#(y) \colon y\in L^1_\F(\mathcal{E}),\:y\leq 0,\:\E_\PP[y|\F]=-1\}\quad \text{ for all }x\in\mathscr{X}.
\] 
In this case, the following are equivalent:

\begin{enumerate}
\item $\rho$ attains its representation;
	\item $\rho$ has the Lebesgue property;
	\item $\rho^\#$ is  stably inf-compact w.r.t. $\sigma_s(L^1_\F(\EE),L^\infty_\F(\EE))$.
\end{enumerate}
\end{thm}

Suppose that $\mathcal{X}=L^p(\mathcal{E})$ with $1\leq p<\infty$.  
In this case, every convex risk measure has the Lebesgue property, is representable and the representation is attained for every $x\in L^p(\mathcal{E})$ (see eg~\cite[Theorem 2.11]{kaina2009convex}). 
Thus, we have:

\begin{thm}\label{thm: repILp}
Suppose that $(p,q)$ are H\"{o}lder conjugates with $1\leq p<\infty$. 
If  $\rho:L^p_\F(\mathcal{E})\rightarrow L^0(\F)$ a conditional risk measure, then $\rho$ has the Lebesgue property, and for every $x\in L^p_\F(\EE)$ there exists $ y\in L^q_\F(\EE)$ with $y\leq 0$ and $\E_\PP[y|\F]=-1$ such that 
$
\rho(x)=\E_\PP[x y|\F]-\rho^\#(y).
$
\end{thm}

If $\mathcal{X}:=L^\phi(\EE)$ with $\phi$ a Young function, due Theorem \ref{prop: transferI} we have that $\mathcal{X}^\#:=L^\phi(\EE)$ and applying Theorem \ref{thm: rep} we have the following:

\begin{thm}\label{thm: repOrlicz}
Let $(\phi,\psi)$ be Young conjugate functions and $\rho:L^\phi_\F(\EE)\rightarrow L^0(\F)$ a conditional risk measure. 
Then $\rho$ is representable, i.e. 
\[
\rho(x)=\esssup\{\E_\PP[x y|\F]-\rho^\#(y)\colon y\in L^\psi_\F(\EE),\:y\leq 0,\E_\PP[y|\F]=-1\}\quad\text{ for all }x\in\mathscr{X}
\]
if and only if $\rho$ is lower semicontinuous w.r.t. $\sigma_s(L^\phi_\F(\EE),L^\psi_\F(\EE))$. 

In that case, the following are equivalent:

\begin{enumerate}
\item $\rho$ attains its representation;
	\item $\rho$ has the Lebesgue property;
	\item $\rho^\#$ is stably inf-compact w.r.t. $\sigma_s(L^\psi_\F(\EE),L^\phi_\F(\EE))$.
\end{enumerate}
\end{thm}

If $\mathcal{X}:=H^\phi(\EE)$ with $\phi$ finite-valued, then  every convex risk measure on $H^\phi(\EE)$ has the Lebesgue property, is representable and the representation is attained for every $x\in H^\phi(\mathcal{E})$ (see eg~\cite[Theorem 4.4]{cheridito02}). 
Thus, we have the following:

 \begin{thm}\label{thm: repOrliczHeart}
Let $(\phi,\psi)$ be Young conjugate functions with $\phi$ finite-valued and    $\rho:H^\phi_\F(\mathcal{E})\rightarrow L^0(\F)$  a conditional risk measure. 
Then, $\rho$ has the Lebesgue property, and for every $x\in H^\phi_\F(\mathcal{E})$ there exists $ y\in L^\psi_\F(\mathcal{E})$ with $y\leq 0$ and $\E_\PP[y|\F]=-1$ such that 
$
\rho(x)=\E_\PP[x y|\F]-\rho^\#(y).
$
\end{thm}

Note that all these theorems are just some examples: we can state a version
of any theorem $\varphi$ on duality theory of convex risk measures and it immediately renders a version for
 conditional risk measures of the form $\llbracket\varphi \rrbracket=I$. 
 
 Also, we would like to point out that the relations in Theorem \ref{prop: transferI} can be easily increased. Moreover, it is also possible to cover more general cases: conditional risk measures with values in $L^0(\F,\mathbb{R}\cup\{+\infty\})$, quasi-convex conditional risk measures and so on. 

\section{Sketch of the proof of the main result}

For saving space, we will give only a sketch of the proof of Theorem~\ref{prop: transferI}. 
A more detailed exposition can be found in \cite[Chapter 4]{zapatathesis}.  

The set of real numbers is a definable notion of ZFC. 
We will denote by $\R_{\A}$ the unique name in $\VAA$ that satisfies the definition of real numbers, which exists due to the transfer and maximum principles. 
Likewise, we will denote by $\N_{\A}$ the unique name in $\VAA$ which satisfies the definition of natural numbers. 

It is well-known that $L^0(\F)$ is a Boolean valued interpretation of the real numbers, see~\cite[Chapter 2, Section 2]{takeuti1978}. 
More precisely, we can state this fact as follows: there is a bijection 
\[
\begin{array}{cccc}
\imath: & L^0(\F) &\longrightarrow & \R_\A{\downarrow}\\
 &\eta &\longmapsto &\eta^\bullet\\
        &u^\circ & \longmapsfrom &u
\end{array}
\]
such that the following is satisfied:

\begin{itemize}
\item[(i)] $\imath(L^0(\F,\N))=\N_{\A}{\downarrow}$ and $(\sum 1_{a_k}n_k)^\bullet=\sum \check{n}_k a_k$ whenever $(n_k)\subset \mathbb{N}$ and $(a_k)\in p(I)$;\footnote{As usual, $L^0(\F,\N)$ denotes the set of classes of equivalence of $\mathbb{N}$-valued measurable functions.}
\item[(ii)] $\llbracket 0^\bullet=0\rrbracket=I$, $\llbracket1^\bullet=1\rrbracket=I$, $\llbracket\eta^\bullet+\xi^\bullet=(\eta+\xi)^\bullet\rrbracket=I$ and $\llbracket\eta^\bullet\xi^\bullet=(\eta\xi)^\bullet\rrbracket=I$ for all $\eta,\xi\in L^0(\F)$; 
\item[(iii)] $\llbracket \eta^\bullet=\xi^\bullet\rrbracket=\bigvee\{a\in\A  \colon 1_a\eta= 1_a\xi\}$ and $\llbracket \eta^\bullet\leq \xi^\bullet\rrbracket=\bigvee\{ a\in\A \colon  1_a\eta\leq 1_a\xi\}$ for all $\eta,\xi\in L^0(\F)$;
\item[(iv)] $\left(\sum 1_{a_k}\eta_k\right)^\bullet=\sum \eta_k^\bullet a_k$ for each $(\eta_k)\subset L^0(\F)$ and $(a_k)\in p(I)$.
\end{itemize}

Now, write $\overline{\R}_{\A}$ for the unique name in $\VAA$ that satisfies the definition of the extended real numbers. 
Using the same techniques as in \cite{takeuti1978} it can be proved that the function $\imath$ extends to a bijection
$$\bar{\imath}:\overline{L^0}(\F)\rightarrow \overline{\R}_{\A}{\downarrow},\quad\eta\mapsto\eta^\bullet$$ such that
\[
\llbracket \eta^\bullet=\xi^\bullet\rrbracket=\bigvee\{a\in\A  \colon 1_a\eta= 1_a\xi\},
\]
where, by convention, we take above $0\cdot(\pm\infty)=0$.

Suppose that $(\eta_n)$ is a sequence in $L^0(\F)$. 
Then, we define $(\eta_\n)_{\n\in L^0(\F,\N)}$ where $\eta_{\n}=\sum_{k\in\N} 1_{\{\n=k\}}\eta_k$.
Then the function 
$$\N_{\A}{\downarrow}\rightarrow \R_{\A}{\downarrow},\quad \n^\bullet\mapsto \eta_{\n}^\bullet,$$ is extensional. 
Due to Proposition \ref{thm: extensional}, we can find a name $v$ with $\llbracket v:\mathbb{N}\to \mathbb{R}\rrbracket=I$  such that $$\llbracket \eta_{\n}^\bullet=v(\n^\bullet)\rrbracket=I,$$
for all $\n$. 
Moreover, we have the following: 

\begin{prop}\cite[Proposition 2.2.1]{takeuti2015two}\label{lem: seq}
If $(\eta_n)$ is a sequence in $L^0(\F)$, then 

$\llbracket(\liminf_n \eta_n)^\bullet=\liminf_n \eta_{n^\circ}^\bullet\rrbracket=I$ and $\llbracket(\limsup_n \eta_n)^\bullet=\limsup_n \eta_{n^\circ}^\bullet\rrbracket=I$. 
In particular, 
$\lim_n \eta_n=\eta$ a.s. if and only if $\llbracket\lim_n \eta^\bullet_{n^\circ}=\eta^\bullet\rrbracket=I$.
\end{prop} 

Suppose that $(X,\Sigma,Q)_{\A}$ is a name for a probability space.\footnote{That $(X,\Sigma,Q)_{\A}$ is a name for a probability space means that $X$ is a name for a set, $\Sigma$ is a name for a $\sigma$-algebra on $X$, $Q$ is a name for a probability measure on $\Sigma$, and $(X,\Sigma,Q)_{\A}$ denotes the corresponding ordered triple within $\VAA$. }   
We can consider the name $L^1(\Sigma)_{\A}$ for the space of classes of equivalence of random variables with finite expectation. 

Gordon \cite[Theorem 5]{gordonCondExp} proved that the conditional expectation $\E_\PP[\cdot|\F]$ from $L^1(\mathcal{E})$ to $L^1(\F)$ is a Boolean valued interpretation of the name for the expectation $\E_Q[\cdot]$ for some probability measure $Q$ within $\VA$. 
 We state this fact in the following proposition, whose self-contained proof can be found in \cite[Section 4.1]{zapatathesis}.

\begin{prop}\label{prop: L1new}
There exists a name $(X,\Sigma,Q)_{\A}$ for a probability space and a bijection
\[
\begin{array}{cccc}
\jmath: & L^1_\F(\mathcal{E}) &\longrightarrow & L^1(\Sigma)_{\A}{\downarrow}\\
 &x &\longmapsto &x^\bullet\\
        &u^\circ & \longmapsfrom &u
\end{array}
\]
such that:
\begin{enumerate}
\item $\jmath$ extends the canonical isomorphism $\imath$;
\item $\llbracket\E_{\PP}[x|\F]^\bullet=\E_Q[x^\bullet]\rrbracket=I$ for all $x\in L^1_\F(\mathcal{E})$;
\item $\llbracket x^\bullet=y^\bullet\rrbracket=\bigvee\{a\in\A\colon 1_a x=1_a y\}$ for all $x,y\in L^1_\F(\mathcal{E})$;
\item $\llbracket x^\bullet\leq y^\bullet\rrbracket=\bigvee\{a\in\A\colon 1_a x\leq 1_a y\}$ for all $x,y\in L^1_\F(\mathcal{E})$;
\item $\llbracket x^\bullet+y^\bullet\rrbracket=\llbracket (x+y)^\bullet\rrbracket$ for all $x,y\in L^1_\F(\mathcal{E})$;

\item $\left(\sum 1_{a_k}x_k\right)^\bullet=\sum x_k^\bullet a_k$ for all $(x_k)\subset L^1_\F(\mathcal{E})$ and $(a_k)\in p(I)$.
\end{enumerate}
\end{prop}

For the forthcoming discussion, we will fix a name for a probability space  $(X,\Sigma,Q)_{\A}$ as in the theorem above. 

Suppose that $S$ is a stable subset of $L^1_\F(\mathcal{E})$. 
Let $S{\uparrow}$ denote the unique representative in $\VAA$ of the name given by the function
\[
\{x^\bullet\colon x\in S\}\longrightarrow\A,\quad x^\bullet\mapsto I.
\]

Using the mixing principle, it is not difficult to prove the following: 
 
\begin{prop}
If $S$ is a stable subset of $L^1_\F(\mathcal{E})$, then $S{\uparrow}$ is a name for a non-empty subset of $L^1_\F(\mathcal{E}){\uparrow}$, and the  map $x\mapsto x^\bullet$ is a bijection from $S$ to $S{\uparrow}{\downarrow}$. 
In particular, we have that $\llbracket L^1_\F(\mathcal{E}){\uparrow}=L^1(\Sigma)_{\A}\rrbracket=I$. 
\end{prop}

Both $\mathscr{X}$ and $\mathscr{X}^\#$ are stable subsets of $L^1_\F(\mathcal{E})$. 
Thus, it makes sense to define  $\mathscr{X}{\uparrow}$ and $\mathscr{X}^\#{\uparrow}$, which are the names that we refer to in the statement of Theorem \ref{prop: transferI}.
 Indeed, 
bearing in mind the properties given in Proposition \ref{prop: L1new}, a standard manipulation of Boolean truth values proves the following:

\begin{prop}\label{prop: dualK} 
 $\mathscr{X}{\uparrow}$ and $\mathscr{X}^\#{\uparrow}$ are names for solid subspaces of $L^1_\F(\mathcal{E}){\uparrow}$ with
 $\llbracket\R\subset\mathscr{X}{\uparrow}\rrbracket=I$ and $\llbracket\R\subset\mathscr{X}^\#{\uparrow}\rrbracket=I$.  
 Moreover, $\llbracket\mathscr{X}^\#{\uparrow}=\mathscr{X}{\uparrow}^\#\rrbracket=I$.
 \end{prop}

If $p$ is a real number with $1\leq p\leq \infty$, we have that its canonical inversion, $\check{p}$ in $\VAA$,  satisfies that $\llbracket \check{p}=p\rrbracket=I$. 
Then we can consider the corresponding name, say $L^{p}(\Sigma)_\A$, for a $L^p$ space within $\VA$. 

Given a Young function $\phi$, consider the function $\eta\mapsto \phi(\eta)\colon L^0(\F,[0,\infty))\rightarrow L^0(\F,[0,\infty])$. 
Due to Proposition \ref{thm: extensional}, we have a name $\tilde{\phi}$ for a Young function. 
Then we can consider the corresponding names for an Orlicz space and an Orlicz-heart space within $\VA$, say  $L^{\tilde{\phi}}(\Sigma)_\A$ and $H^{\tilde{\phi}}(\Sigma)_\A$, respectively. 

The following result tells us that $L^p$, Orlicz  and Orlicz-heart type modules can be interpreted as classical $L^p$, Orlicz  and Orlicz-heart spaces within $V^{(\mathcal{A})}$.  
Actually, general $L^p$ type modules $L^p(\Phi)$ where $\Phi$ is a Maharam operator were introduced and their Boolean valued interpretation was provided in \cite[4.2.2]{kusraev1985vector}. 
In fact, a Maharam operator can be viewed as an abstract conditional expectation, see \cite[Sections 5.2--5.4]{kusraev2014boolean}; moreover, if $\Phi$ is the conditional expectation, then $L^p(\Phi)$ is precisely $L^p_\F(\mathcal{E})$.

By manipulation of Boolean truth values, and bearing in mind the properties given in Proposition \ref{prop: L1new},  we can check the following:

\begin{prop}\label{cor: Orlicz}
If $1\leq p\leq \infty$, then $\llbracket L^p_\F(\EE){\auu}=L^{p}(\Sigma)_\A\rrbracket=I$.   
If $\phi$ is a Young function, then $\llbracket L^\phi_\F(\EE){\auu}=L^{\tilde{\phi}}(\Sigma)_\A\rrbracket=I$ and 
$\llbracket H^\phi_\F(\EE){\auu}=H^{\tilde{\phi}}(\Sigma)_\A\rrbracket=I$.
\end{prop}

Suppose that $(x_n)$ is a sequence in $L^1_\F(\mathcal{E})$. 
For each $\mathfrak{n}\in L^0(\F,\N)$ we define $x_{\mathfrak{n}}:=\sum_{k\in\N}1_{\{\n=k\}}x_k$. 
Then the function
\[
\N_{\A}{\downarrow}\longrightarrow L^1_\F(\mathcal{E}){\uparrow}{\downarrow},\quad \n^\bullet\mapsto x_{\n}^\bullet
\]
is extensional. 
Due to Proposition \ref{thm: extensional} we can find a name $(x_n){\uparrow}$ for a sequence in $L^1_\F(\mathcal{E}){\uparrow}$.
In addition, a standard manipulation of Boolean truth values proves the following:

\begin{prop}\label{prop: seqL1}
Let $(x_n)$ be  a sequence in $\mathscr{X}$ such that $|x_n|\leq y$ for some $y\in \mathscr{X}$. 
Then $\llbracket(\liminf_n x_n)^\bullet=\liminf_n x_{n^\circ}^\bullet\rrbracket=I$ and $\llbracket(\limsup_n x_n)^\bullet=\limsup_n x_{n^\circ}^\bullet\rrbracket=I$.
In particular, $\lim_n x_n=x$ a.s. in $\mathscr{X}$ if and only if $\llbracket\lim_n x^\bullet_{n^\circ}=x^\bullet\textnormal{ a.s. in }\mathscr{X}{\auu}\rrbracket=I.$

\end{prop}

A function $f:\mathscr{X}\rightarrow\overline{L^0}(\F)$ is said to have the \emph{local property} if $1_a f(x)=1_a f(1_a x)$ for all $a\in\mathcal{A}$ and $x\in\mathscr{X}$. 
It is not difficult to verify that if $f$ has the local property, then the function
  
$$\mathscr{X}{\uparrow}{\downarrow}\longrightarrow \overline{\R}_{\A}{\downarrow},\quad x^\bullet\mapsto f(x)^\bullet,$$ is extensional. 
Thus, we can find a name $f{\uparrow}$ for a function from $\mathscr{X}{\uparrow}$ to $\overline{\R}_{\A}$ such that 
$\llbracket f{\uparrow}(x^\bullet)=f(x)^\bullet \rrbracket=I$ for all $x\in\mathscr{X}$.

The following is a consequence of Propositions \ref{lem: seq}:

\begin{prop}\label{prop: FatLeb}
Let $f:\mathscr{X}\rightarrow \overline{L^0}(\F)$ be  a function with the local property. 
Then
\begin{enumerate}
\item $f$ has the Fatou property iff $\llbracket f{\auu}\textnormal{ has the Fatou property}\rrbracket=I$;
\item $f$ has the Lebesgue property iff $\llbracket f{\auu}\textnormal{ has the Lebesgue property}\rrbracket=I$.
\end{enumerate}
\end{prop}

Proposition \ref{thm: forallExists} together with the fact that $x\mapsto x^\bullet$ is a bijection from $S$ to $S{\uparrow}{\downarrow}$ allows us to prove the following:

\begin{prop}\label{prop: stabSet}
Let $S\subset\mathscr{X}$ be stable, let $f:\mathscr{X}\rightarrow \overline{L^0}(\F)$ be a function with the local property. 
Then 
$$\left\llbracket\underset{u\in S{\uparrow}}\sup f{\uparrow}(u)=\left(\underset{x\in S}\esssup f(x)\right)^\bullet\right\rrbracket=I.$$
\end{prop}

Suppose that $\mathscr{B}$ is a stable collection consisting of stable subsets of $L^1_\F(\mathcal{E})$. 

Let $\mathscr{B}{\Uparrow}$ denote the unique name in $\VAA$ equivalent to the name given by the function
\[
\left\{S{\uparrow}  \colon S\in\mathscr{B} \right\}\longrightarrow\A,\quad S{\uparrow}\mapsto I.
\] 

By means of a manipulation of Boolean truth values using the mixing principle, one can prove the following:

\begin{prop}
Let $\mathscr{B}$ be a stable collection consisting of stable subsets of $L^1_\F(\mathcal{E})$. 
Then $\mathscr{B}{\Uparrow}$ is a name for a non-empty collection of non-empty subsets of $L^1_\F(\mathcal{E}){\uparrow}$, and 
the map $S\mapsto S{\uparrow}$ is a bijection from $\mathscr{B}$ to $\mathscr{B}{\Uparrow}{\downarrow}$. 
\end{prop}

Notice that both $\mathscr{B}_{\langle\mathscr{X},\mathscr{X}^\#\rangle}$ and $\mathscr{B}_{\langle\mathscr{X}^\#,\mathscr{X}\rangle}$ are stable collections. 
Then it makes sense to define $\mathscr{B}_{\langle\mathscr{X},\mathscr{X}^\#\rangle}{\Uparrow}$ and $\mathscr{B}_{\langle\mathscr{X}^\#,\mathscr{X}\rangle}{\Uparrow}$. 

Dual systems of modules were introduced and studied in \cite{kusraev1985vector}. 
In addition, their Boolean valued representation can be found in \cite[Theorem 3.3.10]{kusraev1985vector}, which covers the stable weak topologies. 
Actually, we have the following result, which can be also proved by adapting the proof of \cite[Proposition 2.3.20]{zapatathesis}:  

\begin{prop}
$\mathscr{B}_{\langle\mathscr{X},\mathscr{X}^\#\rangle}{\Uparrow}$  (resp. $\mathscr{B}_{\langle\mathscr{X}^\#,\mathscr{X}\rangle}{\Uparrow}$) is a name for a topological base of the weak topology $\sigma(\mathscr{X}{\uparrow},\mathscr{X}^\#{\uparrow})$  (resp. $\sigma(\mathscr{X}{\uparrow}^\#,\mathscr{X}{\uparrow})$) within $\VA$. 
\end{prop}

Next, we deal only with the topology $\sigma_s(\mathscr{X},\mathscr{X}^\#)$, but the following results are also valid for the topology $\sigma_s(\mathscr{X}^\#,\mathscr{X})$. 

The next proposition can be proved by manipulation of the Boolean truth values as in \cite[Proposition 2.3.4]{zapatathesis} and 
 \cite[Corollary 2.3.1]{zapatathesis}.
 
\begin{prop}
Let $S$ be a stable subset of   $\mathscr{X}$.  
Then:
\begin{enumerate}
\item $S$ is open w.r.t. $\sigma_s(\mathscr{X},\mathscr{X}^\#)$ iff 
$\llbracket S{\uparrow}\textnormal{ is open w.r.t. }\sigma(\mathscr{X}{\uparrow},\mathscr{X}^\#{\uparrow})\rrbracket=I$;
\item $S$ is closed w.r.t. $\sigma_s(\mathscr{X},\mathscr{X}^\#)$ iff $\llbracket S{\uparrow}\textnormal{ is closed w.r.t. }\sigma(\mathscr{X}{\uparrow},\mathscr{X}^\#{\uparrow})\rrbracket=I$;
\item $S$ is stably compact  w.r.t. $\sigma_s(\mathscr{X},\mathscr{X}^\#)$ iff $\llbracket S{\uparrow}\textnormal{ is compact w.r.t. }\sigma(\mathscr{X}{\uparrow},\mathscr{X}^\#{\uparrow})\rrbracket=I$.
\end{enumerate} 

\end{prop}  

A function $f:\mathscr{X}\rightarrow \overline{L^0}(\F)$ is said to be \emph{proper} if $f(x)>-\infty$ and there exists $x_0\in \mathscr{X}$ such that $f(x_0)\in L^0(\F)$.

As a consequence of the previous result and by means of a manipulation of the Boolean truth values as in \cite[Proposition 2.3.11]{zapatathesis}, we obtain the following:

\begin{prop}\label{prop:semInfCompact}
Let $f:\mathscr{X}\rightarrow \overline{L^0}(\F)$ be  a proper function with the local property. 
Then:
\begin{enumerate}
\item $f$ is lower semicontinuous  w.r.t. $\sigma_s(\mathscr{X},\mathscr{X}^\#)$ iff 

$\llbracket f\textnormal{ is lower semicontinuous w.r.t. }\sigma(\mathscr{X}{\uparrow},\mathscr{X}^\#{\uparrow})\rrbracket=I$;
\item $f$ is stably inf-compact  w.r.t. $\sigma_s(\mathscr{X},\mathscr{X}^\#)$ iff 

$\llbracket f\textnormal{ is inf-compact w.r.t. }\sigma(\mathscr{X}{\uparrow},\mathscr{X}^\#{\uparrow})\rrbracket=I$.
\end{enumerate}
\end{prop}
  
At this point, we can already prove  Theorem~\ref{prop: transferI}. 
Namely, suppose that $\rho:\mathscr{X}\rightarrow L^0(\F)$ is a conditional risk measure. 
Since $\rho$ is $L^0(\F)$-convex, we know from \cite[Theorem 3.2]{kupper03} that $\rho$ has the  local property. 
Then, we have a name $\rho{\uparrow}$ for a function from $\mathscr{X}{\uparrow}$ to $\R_{\A}$ such that 
$\llbracket\rho{\uparrow}(x^\bullet)=\rho(x)^\bullet \rrbracket=I$ for all $x\in\mathscr{X}$. 
Moreover, it can be computed that $\rho{\uparrow}$ is a name for a convex risk measure.

We also have that $\rho^\#$ is a proper function with the local property; thus, we can find a name $\rho^\#{\uparrow}$ for a proper function from $\mathscr{X}^{\#}$ to $\overline{\R}_{\A}$ so that 
$\llbracket\rho^\#{\uparrow}(y^\bullet)=\rho^\#(y)^\bullet \rrbracket=I$ for all $y\in\mathscr{X}^\#$. 
In addition, as a consequence of Proposition \ref{prop: stabSet}, one has that $\llbracket\rho^\#{\uparrow}=\rho{\uparrow}^\#\rrbracket=I$.

Finally, note that 1 in Theorem ~\ref{prop: transferI} is a consequence of Proposition \ref{prop: stabSet}; 2  in Theorem ~\ref{prop: transferI} is clear from Proposition~\ref{prop: L1new}; 3--4 in Theorem ~\ref{prop: transferI} is precisely Proposition \ref{prop: FatLeb}; 5 in Theorem ~\ref{prop: transferI} is clear from Proposition~\ref{prop: L1new}; 
6--7 in Theorem ~\ref{prop: transferI} is just Proposition \ref{prop:semInfCompact}; and finally we obtain 8--10 in Theorem ~\ref{prop: transferI} from Proposition \ref{cor: Orlicz}.\\

\textbf{Acknowledgment.}  
This paper is based on some parts of the doctoral thesis of the author. 
He would like to thank his mentor José Orihuela his support and  Antonio Avilés López for fruitful discussions.
Likewise, he would like to express his gratitude to Anatoly G. Kusraev and Sem\"{e}n S. Kutateladze for  taking the time to review his thesis, supporting his work and pointing out some extremely nice references in Boolean valued analysis.


\end{document}